%% file: agt-5-3.tex
\documentclass{gtart_h}
\input agtout

\lognumber{3}
\volumenumber{5}
\volumeyear{2005}
\papernumber{3}
\pagenumbers{31}{51}
\received{21 May 2004} 
\revised{23 December 2004}
\accepted{5 January 2005}
\published{7 January 2005}

\usepackage{amsmath,amssymb,amscd}

\newtheorem{prop}[equation]{Proposition}

\newtheorem{thm}[equation]{Theorem}
\newtheorem{cor}[equation]{Corollary}
\newtheorem{lem}[equation]{Lemma}

\theoremstyle{definition}

\newtheorem{defn}[equation]{Definition}
\newtheorem{rem}[equation]{Remark}
\newtheorem{exa}[equation]{Example}

\numberwithin{equation}{section}

\newcommand{\sands}{\mbox{$\quad\text{and}\quad$}}

\newcommand{\Hom}{\operatorname{Hom}}

\newcommand{\Aut}{\operatorname{Aut}}
\newcommand{\letbe}{\mathbin{:\!\raisebox{-.32pt}{=}}}

\newcommand{\EM}{Eilenberg-Mac Lane}

\newcommand{\bkss}{Bousfield-Kan spectral sequence}

\newcommand{\bN}{\mathbb{N}}

\newcommand{\bQ}{\mathbb{Q}}

\newcommand{\bZ}{\mathbb{Z}}

\newcommand{\fcat}[2]{\mbox{\raisebox{1pt}{\rm\scriptsize [}{\sc #1},\hspace{1pt}{\sc #2}\raisebox{1pt}{\rm\scriptsize ]}}\hspace{1pt}}

\newcommand{\llongrightarrow}{\relbar\joinrel\longrightarrow}

\newcommand{\llongleftarrow}{\longleftarrow\joinrel\relbar}

\newcommand{\rla}{\begin{picture}(40,7.5)\put(6,-2.5){$\longleftarrow\joinrel\hspace{-1pt}\relbar\joinrel\hspace{-1pt}\relbar$}\put(6,2.5){$\relbar\joinrel\hspace{-1pt}\relbar\joinrel\hspace{-1pt}\longrightarrow$}
\end{picture}}
\newcommand{\que}{\mathop{\rla}}
\newcommand{\under}{\!\downarrow\!}
\newcommand{\nni}{\not\ni}

\newcommand{\HSM}{\mbox{$\operatorname{End}_{ho}$}}
\newcommand{\HA}{\mbox{$\operatorname{Aut}_{ho}$}}
\newcommand{\col}{\operatorname{colim}}

\newcommand{\hoc}{\operatorname{hocolim}}

\newcommand{\daja}{Davis-Januszkiewicz}
\newcommand{\daaja}{Davis and Januszkiewicz}
\newcommand{\chasc}{Chacholski and Scherer}

\newcommand{\CPI}{\mbox{${\it CP\/}^\infty$}}

\newcommand{\GL}{\mbox{\it GL\/}}

\newcommand{\srf}[1]{\mbox{${\text{\it SR}^F}$}}

\newcommand{\etminu}{\!\setminus\!}

\newcommand{\APL}{\mbox{$A_{P\!L}$}}

\newcommand{\cat}[1]{\mbox{\sc #1}}

\newcommand{\del}{\mbox{\footnotesize{$\Delta$}}}

\newcommand{\ssupseteq}{\mbox{$\scriptscriptstyle \supseteq$}}

\newcommand \lra{\longrightarrow}

\begin{document}

\title[Homotopy types, formality and rationalisation]
{On Davis-Januszkiewicz homotopy types I;\\formality 
and rationalisation}

\authors{Dietrich Notbohm\\Nigel Ray}

\address{Department of Mathematics and Computer Science, University of
Leicester\\University Road, Leicester LE1 7RH, UK}
\secondaddress{Department of Mathematics, University of Manchester\\Oxford
Road, Manchester M13 9PL, UK}
\asciiaddress{Department of Mathematics and Computer Science, University of
Leicester\\University Road, Leicester LE1 7RH, UK\\and\\Department 
of Mathematics, University of Manchester\\Oxford
Road, Manchester M13 9PL, UK}
\asciiemail{dn8@mcs.le.ac.uk, nige@ma.man.ac.uk}
\gtemail{\mailto{dn8@mcs.le.ac.uk}{\rm\qua and\qua}\mailto{nige@ma.man.ac.uk}}

\asciikeywords
{Colimit, formality, Davis-Januszkiewicz space, homotopy colimit, model
category, rationalisation, Stanley-Reisner algebra}

\begin{abstract}
For an arbitrary simplicial complex $K$, Davis and Januszkiewicz have
defined a family of homotopy equivalent CW-complexes whose integral
cohomology rings are isomorphic to the Stanley-Reisner algebra of
$K$. Subsequently, Buchstaber and Panov gave an alternative
construction (here called $c(K))$, which they showed to be homotopy
equivalent to Davis and Januszkiewicz's examples. It is therefore
natural to investigate the extent to which the homotopy type of a
space is determined by having such a cohomology ring. We begin this
study here, in the context of model category theory. In particular, we
extend work of Franz by showing that the singular cochain algebra of
$c(K)$ is formal as a differential graded noncommutative algebra. We
specialise to the rationals by proving the corresponding result for
Sullivan's commutative cochain algebra, and deduce that the
rationalisation of $c(K)$ is unique for a special family of complexes
$K$. In a sequel, we will consider the uniqueness of $c(K)$ at each
prime separately, and apply Sullivan's arithmetic square to produce
global results for this family.
\end{abstract}
\asciiabstract{%
For an arbitrary simplicial complex K, Davis and Januszkiewicz have
defined a family of homotopy equivalent CW-complexes whose integral
cohomology rings are isomorphic to the Stanley-Reisner algebra of
K. Subsequently, Buchstaber and Panov gave an alternative construction
(here called c(K)), which they showed to be homotopy equivalent to
Davis and Januszkiewicz's examples. It is therefore natural to
investigate the extent to which the homotopy type of a space is
determined by having such a cohomology ring. We begin this study here,
in the context of model category theory. In particular, we extend work
of Franz by showing that the singular cochain algebra of c(K) is
formal as a differential graded noncommutative algebra. We specialise
to the rationals by proving the corresponding result for Sullivan's
commutative cochain algebra, and deduce that the rationalisation of
c(K) is unique for a special family of complexes K. In a sequel, we
will consider the uniqueness of c(K) at each prime separately, and
apply Sullivan's arithmetic square to produce global results for this
family.}

\maketitle

\section{Introduction}\label{intro}

Over the last decade, work of \daaja\ \cite{daja:cpc} has popularised
homotopy theoretical aspects of toric geometry amongst algebraic
topologists. The results of \cite{daja:cpc} have been surveyed by
Buchstaber and Panov in \cite{bupa:taa}, where several further
applications were developed. Their constructions have led us to consider
the uniqueness of certain associated homotopy types, and our aim is to
begin that study here; we focus on general issues of formality, and
deduce that the rationalisations of a family of special cases are
unique. In a sequel \cite{nora:djh2}, we discuss the problem prime by
prime, and obtain global results for members of the family by appeal to
Sullivan's arithmetic square. The general problem of uniqueness appears
to be of considerable difficulty.

We work over an arbitrary commutative ring $R$ with identity, and
consider a universal set $V$ of {\it vertices} $v_1$, \dots, $v_m$,
ordered by their subscripts. The vertices masquerade as
algebraically independent variables, which generate a graded
polynomial algebra $S_R(V)$ over $R$. The grading is defined by
assigning each of the generators a common dimension, which we
usually take to be $2$. A function $M\colon V\rightarrow \bN$ is
known as a {\it multiset\/} on $V$, with cardinality
$|M|\letbe\sum_jM(v_j)$; it may be represented by the monomial
$v_M\letbe\prod_V\smash{v^{M(v)}}$, or by the $n$-tuple of
constituent vertices $(v_{j_1},\dots,v_{j_n})$, where
$j_1\leq\dots\leq j_n$ and $n=|M|$. So $S_R(V)$ is generated
additively by the $v_M$, and $v_M$ is squarefree precisely when $M$
is a genuine subset.

A simplicial complex $K$ on $V$ consists of a finite set of faces
$\sigma\subseteq V$, closed with respect to the formation of subsets
$\rho\subseteq\sigma$. Alternatively, we may interpret $K$ as the set of
squarefree monomials $v_\sigma\letbe\prod_\sigma v$, which is closed
under factorisation. Every simplicial complex generates a simplicial set
$K_\bullet$; for each $n\geq 0$, the $n$-simplices $K_n$ contain all $M$
of cardinality $n+1$ whose support is a face of $K$. The face and
degeneracy operators delete and repeat the appropriate vertices
respectively.

The Stanley-Reisner algebra $R[K]$, otherwise known as the {\it face
ring\/} of $K$, is an important combinatorial invariant. It is defined
as the quotient
\begin{equation}\label{defsra}
S_R(V)/(v_U:U\notin K),
\end{equation}
and is therefore generated additively by the simplices of $K_\bullet$.
The algebraic properties of $R[K]$ encode a host of combinatorial
features of $K$, and are discussed in detail by Bruns and Herzog
\cite{brhe:cmr} and Stanley \cite{st:cca}, for example. If $K$ is the
simplex on $V$, then $R[K]$ is the polynomial algebra $S_R(V)$.

For each $K$, \daaja\ defined the notion of a toric space over the cone
on the barycentric subdivision of $K$, and showed that the cohomology of
such a space is related to $R[K]$. The relationship follows from their
application of the Borel construction, which creates a family of spaces
whose cohomology ring (with coefficients in $R$) is isomorphic to the
Stanley-Reisner algebra. All members of the family are homotopy
equivalent to a certain universal example, and we refer to any space
which shares their common homotopy type as a \daja\ space. The
isomorphisms equip $R[K]$ with a natural grading, which agrees with that
induced from $S_R(V)$.  Subsequently, Buchstaber and Panov
\cite{bupa:taa} defined a CW-complex whose cohomology ring is also
isomorphic to $R[K]$. They confirmed that their complex is a \daja\
space by giving an explicit homotopy equivalence with the universal
example. In \cite{paravo:csa}, their space is described as the pointed
colimit $\col^+B^K$ of a certain $\cat{cat}(K)$-diagram $B^K$, which
assigns the cartesian product $B^\sigma$ to each face $\sigma$ of
$K$. Here $B$ denotes the classifying space of the circle, or $\CPI$,
and $\cat{cat}(K)$ is the category of faces and inclusions.

We say that a space $X$ {\it realises\/} the Stanley-Reisner algebra of
$K$ whenever there is an algebra isomorphism $H^*(X;R)\cong R[K]$. We
denote the rationalisation of $X$ by $X_0$, and write $\HA(X)<\HSM(X)$
for the homotopy classes of self-equivalences of $X$, considered as a
subgroup of the homotopy classes of self-maps with respect to
composition. 

The contents of each section are as follows.

In Section \ref{ba} we describe our notation and prerequisites,
including those aspects of model category theory which provide a useful
context for exponential diagrams and their cohomology. We also explain
why it is equally acceptable to work with the unpointed colimit
$c(K)\letbe\col B^K$. We introduce the Stanley-Reisner algebra in
Section \ref{incoli}, and show that the \bkss\ for $H^*(c(K);R)$
collapses by analysing higher limits of certain
$\cat{cat}(K)$-diagrams. In Section \ref{info} we apply similar
techniques to prove the formality of the singular cochain algebra
$C^*(c(K);R)$. Finally, we specialise to the case $R=\bQ$ in Section
\ref{rafo}, where we confirm that Sullivan's commutative cochain algebra
$\APL(c(K))$ is formal in the commutative sense. We deduce that $\bQ[K]$
determines the rationalisation $c(K)_0$ whenever it is a complete
intersection, and discuss the corresponding automorphism group
$\HA(c(K)_0)$.

The authors are especially grateful to the organisers of the
International Conference on Algebraic Topology, which was held on the
Island of Skye in June 2001. The Conference provided the opportunity
for valuable discussion with several colleagues, amongst whom Octavian
Cornea, Kathryn Hess, and Taras Panov deserve special mention. Without
that remarkable and stimulating environment, our work would not have
begun. We are more recently indebted to John Greenlees, who brought to
our attention a fundamental misconception in a previous version of
this article, and to the referee, for further improvements. We also
thank the London Mathematical Society, whose support of the
Transpennine Topology Triangle has enabled our collaboration to
continue and develop.

\section{Background}\label{ba}

We begin by establishing our notation and prerequisites, recalling
various aspects of \daja\ spaces. We refine results of
\cite{paravo:csa} in the context of model category theory, referring
readers to \cite{dwsp:htm} and \cite{ho:mc} for background details.
Following \cite{vo:cct}, we adopt the model category $\cat{top}$ of
$k$-spaces and continuous functions as our topological workplace. Weak
equivalences induce isomorphisms in homotopy, fibrations are Serre
fibrations, and cofibrations have the left lifting property with
respect to acyclic fibrations. Every function space $Y^X$ is endowed
with the corresponding $k$-topology. Many of the spaces we consider
have a distinguished basepoint $*$, and we write $\cat{top}_+$ for the
model category of pairs $(X,*)$ and basepoint preserving maps.  We
usually insist that the inclusion $*\rightarrow X$ be a cofibration,
in which case $X$ is {\it well-pointed}; this is automatic when $X$ is
a CW-complex and $*$ its $0$--skeleton.

Given a small category $\cat{a}$, we refer to a covariant functor
$D\colon\cat{a}\rightarrow\cat{r}$ as an {\it $\cat{a}$-diagram} in
$\cat{r}$. Such diagrams are themselves the objects of a category
$\fcat{a}{r}$, whose morphisms are natural transformations of
functors. We may interpret any object $X$ of $\cat{r}$ as a constant
diagram, which maps every object of $\cat{a}$ to $X$ and every morphism
to the identity. 
\begin{exa}\label{simco}
{\rm (1)}\qua For each integer $n\geq 0$, the category $\cat{ord}(n)$ has
objects $0$, $1$, \dots, $n$, equipped with a single morphism
$k\rightarrow m$ when $k\leq m$. An $\cat{ord}(n)$-diagram
\begin{equation}\label{ordndiags}
X_0\stackrel{f_1}{\longrightarrow}X_1\stackrel{f_2}{\longrightarrow}
\dots \stackrel{f_n}{\longrightarrow}X_n
\end{equation}
consists of $n$ composable morphisms in $\cat{r}$.

{\rm(2)}\qua The category $\del$ has objects $(n)\letbe\{0,1,\dots,n\}$ for $n\geq
0$, and morphisms the nondecreasing functions; then $\del^{op}$- and
$\del$-diagrams are simplicial and cosimplicial objects of $\cat{r}$
respectively. In particular, $\varDelta\colon\del\rightarrow\cat{top}$
is the cosimplicial space which assigns the standard $n$-simplex
$\varDelta(n)$ to each object $(n)$.
\end{exa}

Given objects $X_0$ and $X_1$ of $\cat{r}$, we write the set of
morphisms $X_0\rightarrow X_1$ as $\cat{r}(X_0,X_1)$; when $\cat{r}$ is
small, the diagrams \eqref{ordndiags} also form a set for every
$n>1$. For any $\cat{r}$ it is often convenient to abbreviate
$\fcat{$\del^{op}$}{r}$ to $\cat{sr}$, and write a generic simplicial
object as $D_\bullet$. In particular, $\cat{sset}$ denotes the category
of simplicial sets $Y_\bullet$.

 From this point on we work with an abstract simplicial complex $K$,
whose faces $\sigma$ are subsets of the vertices $V$. We assume that
the empty face belongs to $K$, and write $K^\times$ when it is
expressly omitted. The integer $|\sigma|-1$ is known as the {\it
dimension} of $\sigma$, and written $\dim\sigma$; its maximum value
$\dim K$ is the dimension of $K$. When $K$ contains every subset of
$V$, we may call it the {\it simplex} $\varDelta(V)$ on $V$. Each face
of $K$ therefore determines a subsimplex $\varDelta(\sigma)$, whose
{\it boundary\/} $\partial(\sigma)$ is the complex obtained by
deleting the subset $\sigma$. We also require the {\it link}
$\ell_K(\sigma)$, whose faces consist of those $\tau\setminus\sigma$
for which $\sigma\subseteq\tau$ in $K$.

\begin{defn}\label{defcatk}
For any simplicial complex $K$, the small category ${\cat{cat}(K)}$
has objects the faces of $K$ and morphisms the inclusions
$i_{\sigma,\tau}\colon\sigma\subseteq\tau$. The empty face
$\varnothing$ is an initial object, and the {\/\em maximal faces}
$\mu$ admit only identity morphisms. The opposite category
$\cat{cat}^{op}(K)$ has morphisms $p_{\tau,\sigma}\letbe
i_{\sigma,\tau}^{op}\colon\tau\supseteq\sigma$, and $\varnothing$ is
final.
\end{defn}
The nondegerate simplices of the nerve $N_\bullet\cat{cat}(K)$ form
the cone on the bary\-centric subdivision $K'$, and those of
$N_\bullet\cat{cat}(K^\times)$ correspond to the subcomplex $K'$. So the
classifying space $B\cat{cat}(K)$, formed by realising the nerve, is a
contractible CW-complex, and $B\cat{cat}(K^\times)$ is a subcomplex
homeomorphic to $K$. We shall study $\cat{cat}(K)$- and
$\cat{cat}^{op}(K)$-diagrams $D$ in various algebraic and topological
categories $\cat{r}$. Usually, $\cat{r}$ is {\it pointed\/} by an object
$*$, which is both initial and final; unless stated otherwise, we then
assume that $D(\varnothing)=*$.

For each face $\sigma$, the {\it overcategories}
$\cat{cat}(K)\under\sigma$ and $\cat{cat}(K)\!\Downarrow\!\sigma$ are
given by restricting attention to those objects $\rho$ for which
$\rho\subseteq\sigma$ and $\rho\subset\sigma$ respectively. The {\it
undercategories} $\sigma\under\cat{cat}(K)$ and
$\sigma\!\Downarrow\!\cat{cat}(K)$ are defined by the objects
$\sigma\subseteq\tau$ and $\sigma\subset\tau$. It follows from the
definitions that
\[
\begin{split}
\cat{cat}(K)\under\sigma&=\cat{cat}(\varDelta(\sigma)),\quad
\cat{cat}(K)\!\Downarrow\!\sigma=\cat{cat}(\partial(\sigma)),\\
&\sigma\under\cat{cat}(K)=\cat{cat}(\ell_K(\sigma))\sands
\sigma\!\Downarrow\!\cat{cat}(K)=\cat{cat}(\ell_K(\sigma)^\times).
\end{split}
\]
Dimension may be interpreted as a functor
$\dim\colon\cat{cat}(K)\rightarrow\cat{ord}(m-1)$, which is a {\it
linear extension\/} in the sense of \cite{ho:mc}; thus $\cat{cat}(K)$ is
{\it direct} and $\cat{cat}^{op}(K)$ is {\it inverse}.

For any model category $\cat{r}$, we may therefore follow Hovey and
impose an associated model structure on the category of diagrams
$\fcat{cat$(K)$}{r}$. Weak equivalences $e\colon C\rightarrow D$ are
given {\it objectwise}, in the sense that $e(\sigma)\colon
C(\sigma)\rightarrow D(\sigma)$ is a weak equivalence in $\cat{r}$ for
every face $\sigma$ of $K$. Fibrations are also given objectwise. To
describe the cofibrations, we consider the restrictions of $C$ and $D$
to the overcategories $\cat{cat}(\partial(\sigma))$, and write $L_\sigma
C$ and $L_\sigma D$ for their respective colimits; $L_\sigma$ is the
{\it latching functor\/} of \cite{ho:mc}. Then $g\colon C\rightarrow D$
is a cofibration whenever the induced maps
\begin{equation}\label{diagcof}
C(\sigma)\amalg_{L_\sigma C}L_\sigma D\longrightarrow D(\sigma)
\end{equation}
are cofibrations in $\cat{r}$ for every face $\sigma$. Alternatively,
the methods of \chasc\ \cite{chsc:htd} lead to the same model structure
on $\fcat{cat$(K)$}{r}$.

There is a dual model category structure on $\fcat{cat$^{op}(K)$}{r}$,
where weak equivalences and cofibrations are given objectwise. To
describe fibrations $f\colon C\rightarrow D$, we consider the
restrictions of $C$ and $D$ to the undercategories
$\cat{cat}^{op}(\partial(\sigma))$, and write $M_\sigma C$ and
$M_\sigma D$ for their respective limits; $M_\sigma$ is the {\it
matching functor\/} of \cite{ho:mc}. Then $f$ is a fibration whenever
the induced maps
\begin{equation}\label{diagf}
C(\sigma)\longrightarrow D(\sigma)\times_{M_\sigma D}M_\sigma C
\end{equation}
are fibrations in $\cat{r}$ for every face $\sigma$.

\begin{defn}\label{defedp}
For any CW-pair $(X,*)$, the {\/\em exponential pair\/} of diagrams
$(X^K,X_K)$ consists of functors
\[
X^K\colon\cat{cat}(K)\longrightarrow\cat{top}_+\sands
X_K\colon\cat{cat}^{op}(K)\longrightarrow\cat{top}_+,
\]
which assign the cartesian product $X^\sigma$ to each face $\sigma$
of $K$. The value of $X^K$ on $i_{\sigma,\tau}$ is the cofibration
$X^\sigma\rightarrow X^\tau$, where the superfluous coordinates are
set to $*$, and the value of $X_K$ on $p_{\tau,\sigma}$ is the
fibration $X^\tau\rightarrow X^\sigma$, defined by projection. The
pair are {\/\em twins}, in the sense that $X_K(p')\cdot
X^K(i)=X^K(j')\cdot X_K(q)$ for every pullback square
\[
\begin{CD}
\sigma\cap\sigma'\!@>j'>>\sigma'\\@V\smash{j}VV@VV\smash{i'}V\\
\sigma @>>i>\tau
\end{CD}
\]
in $\cat{cat}(K)$, where $p'=(i')^{op}$ and $q=j^{op}$.
\end{defn}
The properties of twin diagrams are analogous to those of a Mackey
functor \cite{gr:art}. They include, for example, the fact that each
$X^K(i)$ has left inverse $X_K(p)$, where $p=i^{op}$. Our applications
in Theorem \ref{hilimze} are reminiscent of \cite{jamc:hdc}, where the
acyclicity of certain Mackey functors is established.

The colimit $\col X^K$ is a subcomplex of $X^V$, whose inclusion $r$
is induced by interpreting the elements $\sigma$ of $K$ as faces of the
$(m-1)$-simplex $\varDelta(V)$. Composing $r$ with any of the natural
maps $X^\sigma\rightarrow\col X^K$ yields the standard inclusion
$X^\sigma\rightarrow X^V$. We note that $\col X^K$ is pointed by
$X^\varnothing$, otherwise known as the basepoint $*$, and is
homeomorphic to the pointed colimit $\col^+X^K$ of \cite{paravo:csa}.

We wish to study homotopy theoretic properties of $\col X^K$ in
favourable cases. Yet the colimit functor behaves particularly
poorly in this context, because objectwise equivalent diagrams may
well have homotopy inequivalent colimits. The standard procedure for
dealing with this situation is to introduce the left derived
functor, known as the {\it homotopy colimit}. Following
\cite{hovo:mts}, for example, $\hoc X^K$ may be described by the
two-sided bar construction $B(*,\cat{cat}(K),X^K)$ in $\cat{top}$.
We note that $\hoc X^K$ is also pointed, and is related to the
pointed homotopy colimit $\hoc^+X^K$ of \cite{paravo:csa} by the
cofibre sequence
\[
B\cat{cat}(K)\longrightarrow\hoc
X^K\stackrel{f}{\longrightarrow}\hoc^+X^K
\]
of \cite{boka:hlc}. Since $B\cat{cat}(K)$ is contractible, $f$ is a weak
equivalence. We may therefore concentrate on $\hoc X^K$, and so avoid
basepoint complications when working with function spaces in
\cite{nora:djh2}.

\begin{lem}\label{expcoft}
Every exponential diagram $X^K$ is cofibrant in $\,\fcat{cat$(K)$}{top}$.
\end{lem}
\begin{proof}
The initial $\cat{cat}(K)$-diagram in $\cat{top}$ is the constant
diagram $*$, so $X^K$ is cofibrant whenever the inclusion $*\rightarrow
X^K$ is a cofibration. By \eqref{diagcof}, it suffices to show that the
map $X^{\partial(\sigma)}\rightarrow X^\sigma$ is a cofibration for
every face $\sigma$ of $K$. But this map includes the fat wedge in the
cartesian product, and the result follows. 
\end{proof}
An immediate consequence of Lemma \ref{expcoft} is that the natural
projection
\begin{equation}\label{proj}
\hoc X^K\longrightarrow\col X^K
\end{equation}
is a homotopy equivalence. This exemplifies one of the fundamental
properties of the homotopy colimit functor, and is sometimes called the
Projection Lemma \cite{wezizi:hcc}.

\section{Integral cohomology and limits}\label{incoli}

In this section we work in the category $\cat{mod}_R$ of $R$-modules,
and study the cohomology of limits of exponential diagrams $B^K$,
where $B$ is the classifying space of the circle. For this case only,
we abbreviate \eqref{proj} to $hc(K)\rightarrow c(K)$. We focus on the
relationship between the Stanley-Reisner algebra $R[K]$ and the \bkss\
for $H^*(hc(K);R)$.

We begin by investigating the cohomology of $c(K)$. To simplify
applications in later sections, we consider an arbitrary pair of
twin diagrams $(D_K,D^K)$,
\begin{equation}\label{dkdk}
D_K\colon\cat{cat}(K)\longrightarrow\cat{mod}_R\sands
D^K\colon\cat{cat}^{op}(K)\longrightarrow\cat{mod}_R.
\end{equation}
Thus $D^K(p')\cdot D_K(i)=D_K(j)\cdot D^K(q')$ for every pullback
square $i\cdot j=i'\cdot j'$ in $\cat{cat}(K)$. In particular, $D^K(p)$
has right inverse $D_K(i)$ for every morphism $p=i^{op}$. Such pairs
arise, for example, from any contravariant functor
$D\colon\cat{top}\rightarrow\cat{mod}_R$, by composing with the
exponential twins of Definition \ref{defedp}. So $(D_K,D^K)=(D\cdot
B_K,D\cdot B^K)$, and functoriality ensures the diagrams are twins. In
this case we may apply $D$ to the natural maps $B^\sigma\rightarrow
c(K)\stackrel{r}{\longrightarrow}B^V$, and obtain homomorphisms
\begin{equation}\label{projs}
D(B^V)\stackrel{D(r)}{\llongrightarrow}D(c(K))
\stackrel{h}{\longrightarrow}\lim D^K
\end{equation}
in $\cat{mod}_R$.

By way of example, we consider the case $D=H^{2j}(-,R)$, for any $j\geq
0$. For every face $\sigma$ of $K$, the space $B^\sigma$ is an \EM\
space $H(\bZ^\sigma,2)$, and may be expressed as the realisation of a
simplicial abelian group $H_\bullet(\bZ^\sigma,2)$ whenever convenient
\cite{ma:soa}. As a CW-complex, the cells of $B^\sigma$ are concentrated
in even dimensions, and correspond to the simplices $v_M$ of
$\varDelta(\sigma)_\bullet$. The cellular cohomology group
$H^{2j}(B^\sigma;R)$ is therefore isomorphic to the free $R$-module
generated by those $v_M$ for which $|M|=j$ and the support of $M$ is a
subset of $\sigma$. The diagram $D^K$ of \eqref{dkdk} becomes
\begin{equation}\label{jcobk}
H^{2j}(B^K;R)\colon\cat{cat}^{op}(K)\longrightarrow\cat{mod}_R,
\end{equation}
whose value on $p_{\tau,\sigma}$ is the homomorphism which fixes $v_M$
whenever the support of $M$ lies in $\sigma$, and annihilates it
otherwise; the right inverse is the inclusion induced by $D_K$. When
$D=H^{2j+1}(-,R)$, the diagram is zero.

In the case of cohomology, we may combine the diagrams \eqref{jcobk}
into a graded version
\begin{equation}\label{grcobk}
H^*(B^K;R)\colon\cat{cat}^{op}(K)\longrightarrow\cat{gmod}_R,
\end{equation}
taking values in the category of graded $R$-modules. The cup product
on each of the constituent submodules $H^*(B^\sigma;R)$ is given by
the product of monomials $x_Lx_M=x_{L+M}$, as follows from the case of
a single vertex. In other words, the cohomology ring $H^*(B^\sigma;R)$
is isomorphic to the polynomial algebra $S_R(\sigma)$. So $H^*(B^K;R)$
actually takes values in the category $\cat{gca}_R$ of graded
commutative $R$-algebras, and maps the morphism $p_{\tau,\sigma}$ to
the projection $S_R(\tau)\rightarrow S_R(\sigma)$; the right inverse
is again inclusion.

The homomorphisms \eqref{projs} may similarly be combined as
\begin{equation}\label{hrprojs}
S_R(V)\stackrel{r^*}{\llongrightarrow}H^*(c(K);R)
\stackrel{h}{\llongrightarrow}\lim H^*(B^K;R),
\end{equation}
where the limit is taken in $\cat{gmod}_R$. Since \eqref{grcobk}
is a diagram of algebras, the limit inherits a multiplicative
structure, and it is equally appropriate to interpret
\eqref{hrprojs} in $\cat{gca}_R$. The composition $h\cdot r^*$
is induced by the projections $S_R(V)\rightarrow S_R(\sigma)$. In this 
case, we make one further observation.
\begin{prop}\label{epi}
The homomorphism $r^*$ is epic, and its kernel is the ideal
$(v_U:U\notin K)$.
\end{prop}
\begin{proof}
In each dimension $2j$, the cells $v_M$ of $B^V$ correspond to the
multisets on $V$ with $|M|=j$. The cells of $c(K)$ form a subset, given
by those $M$ whose support is a face of $K$. Hence $r^*$ is epic, and
its kernel is generated by the remaining cells. These coincide with the 
$2j$--dimensional additive generators of the ideal $(v_U:U\notin K)$.
\end{proof}
So there is an isomorphism $R[K]\cong H^*(c(K);R)$ of the
Stanley-Reisner algebra \eqref{defsra}, which plays a central r\^ole in
\cite{bupa:taa}.

Returning to our study of the twins $(D_K,D^K)$, the following
definition identifies a further important property.
\begin{defn}\label{deffat}
A diagram $F^K\colon\cat{cat}^{op}(K)\rightarrow\cat{mod}_R$ of
$R$-modules (graded or otherwise) is {\/\em fat} if the natural map
$F^K(\sigma)\rightarrow\lim F^{\partial(\sigma)}$ is an epimorphism for
every face $\sigma$ of $K$.
\end{defn}
The terminology acknowledges the relationship between
$\partial(\sigma)$ and the fat wedge described in Lemma \ref{expcoft}.
\begin{lem}\label{dkfat}
The twin $D^K$ is fat.
\end{lem}
\begin{proof}
We consider an arbitrary face $\rho$ of $K$, whose vertices we label
$w_k$ for $0\leq k\leq d$; thus $d=\dim\rho$. We write
$\mu_k\letbe\rho\setminus w_k$ for the maximal faces of
$\partial(\rho)$, and abbreviate the morphism $p_{\rho,\mu_k}$ to $p_k$
for $0\leq k\leq d$.

The definition of $\lim$ ensures that $L\letbe\lim
D^{\partial(\rho)}$ appears in an exact sequence
\[
0\longrightarrow L\longrightarrow
\prod_{\rho\supset\sigma}D^K(\sigma)\stackrel{\delta}{\longrightarrow}
\prod_{\rho\supset\tau\supset\sigma}D^K(\sigma),
\]
where $\delta(u)(\tau\supset\sigma)=
u(\sigma)-D^K(p_{\tau,\sigma})u(\tau)$ for any
$u\in\prod_{\rho\supset\sigma}D^K(\sigma)$. Hence $u\in L$ is
determined by the values $u(\mu_k)$. The natural projection
$D^K(\rho)\rightarrow\prod_{\rho\supset\sigma}D^K(\sigma)$ therefore
factors through $L$, and it remains for us to find $u(\rho)\in
D^K(\rho)$ such that $D^K(p_k)(u(\rho))=u(\mu_k)$ for every $0\leq
k\leq d$.

We define $u(\rho)\letbe\sum_{\rho\supset\sigma}
(-1)^{|\rho\setminus\sigma|+1}D_K(i_{\sigma,\rho})u(\sigma)$. The
fact that $D_K$ and $D^K$ are twins implies that
\[
D^K(p_k)D_K(i_{\sigma,\rho})u(\sigma)=
\begin{cases}
D_K(i_{\sigma\setminus w_k,\mu_k})u\big(\sigma\etminu w_k\big)&\text{if
$w_k\in\sigma$}\\ 
D_K(i_{\sigma,\mu_k})u(\sigma)&\text{otherwise},
\end{cases}
\]
for every $0\leq k\leq d$; thus $D^K(p_k)u(\rho)$ is given by 
\[
\sum_{\sigma\nni
w_k}(-1)^{|\rho\setminus\sigma|+1}D_K(i_{\sigma,\mu_k})u(\sigma)
+\sum_{\tau\ni
w_k}(-1)^{|\rho\setminus\tau|+1}D_K(i_{\tau\setminus
w_k,\mu_k})u(\rho\etminu w_k).
\]
But we may write $u(\sigma)$ as $u\big((\sigma\cup w_k)\etminu
w_k\big)$ for any $\sigma\nni w_k$ other than $\mu_k$. So the
summands cancel in pairs, leaving $u(\mu_k)$ as required.
\end{proof}

For cohomology, Lemma \ref{dkfat} contributes to our analysis of
$c(K)$. The homotopy equivalence \eqref{proj} provides a {\it cohomology
decomposition\/} \cite{dw:csh}, in the sense that the cohomology algebra
$H^*(c(K);R)$ may be computed by the Bousfield-Kan spectral sequence
\cite{boka:hlc}
\[
E_2^{i,j}\Longrightarrow H^{i+j}(hc(K);R),
\]
where $E_2^{i,j}$ is isomorphic to the $i\/$th derived functor
$\lim^iH^j(B^K;R)$ for every $i,j\geq 0$. The vertical edge homomorphism
coincides with the map $h$ of \eqref{hrprojs}. Lemma \ref{dkfat} is
required for our computation of these limits, and Corollary
\ref{codesharp} will confirm that the cohomology decomposition is {\it
sharp\/} in Dwyer's language. Our proof uses the calculus of functors and
their limits; the appropriate prerequisites may be deduced from Gabriel
and Zisman \cite[Appendix II \S3]{gazi:cfh}, by dualising their results
for colimits.

In particular, we follow \cite{gazi:cfh} (as expounded in \cite{ol:hls},
for example) by calculating $\lim^iD^K$ as the $i\,$th cohomology group
of a certain cochain complex $\big(C^*(D^K),\delta\big)$ of
$R$-modules. The groups are defined by
\[
C^n(D^K)\;\letbe
\prod_{\sigma_0\ssupseteq\dots\ssupseteq\sigma_n}D^K(\sigma_n)
\quad\text{for $n\geq 0$}\,,
\]
and the differential $\delta\letbe\sum_{k=0}^n(-1)^k\delta^k$ is
defined on $u\in C^n(D^K)$ by
\[
\delta^k(u)(\sigma_0\supseteq\dots\supseteq\sigma_{n+1})\letbe
\begin{cases}
u(\sigma_0\supseteq\dots\supseteq\widehat{\sigma}_k
\supseteq\dots\supseteq\sigma_{n+1})&\text{for $k\neq n+1$}\\
D^K(p_{\sigma_n,\sigma_{n+1}})u(\sigma_0\supseteq\dots
\supseteq\sigma_n)&\text{for $k=n+1$}.
\end{cases}
\]
We may replace $C^*(D^K)$ by its quotient $N^*(D^K)$ of normalised
cochains, for which the faces $\sigma_0\supset\dots\supset\sigma_n$ are
required to be distinct. 

\begin{lem}\label{atomic}
Given a maximal face $\mu$ of $K$, and a diagram
$D\colon\cat{cat}^{op}(K)\rightarrow\cat{mod}_R$ such that
$D(\sigma)=0$ for all $\sigma\neq\mu$, then 
\[
{\textstyle\lim^i}D=
\begin{cases}
D(\mu)&\text{for $i=0$}\\
0&\text{for $i>0$}.
\end{cases}
\]
\end{lem}
\begin{proof}
Since $\mu$ is maximal, the only morphism $\sigma\supseteq\mu$ is
the identity. So the normalised chain complex $N^*(D)$ is $D(\mu)$
in dimension $0$, and $0$ in higher dimensions, as required.
\end{proof}

\begin{thm}\label{hilimze}
For any fat diagram $F^K\colon\cat{cat}^{op}(K)\rightarrow\cat{mod}_R$,
we have that $\lim^iF^K=0$ for all $i>0$; in particular, $\lim^iD^K=0$
for every twin $D^K$.
\end{thm}
\begin{proof}
We proceed by induction on the total number of faces $f(K)$; the
result obviously holds for the initial example $K=\varnothing$, where
$f(K)=0$. Our inductive hypothesis is that $\lim^iF^K$ vanishes
whenever $K$ satisfies $f(K)\leq f$.

We therefore consider an arbitrary complex $K$ with $f(K)=f+1$, and
write $J\subset K$ for the subcomplex obtained by deleting a single
maximal face $\mu$. The inclusion of $J$ defines a functor
$G\colon\cat{cat}^{op}(J)\rightarrow\cat{cat}^{op}(K)$, whose
induced functor $G^*\colon\fcat{cat$^{op}(K)$}{mod$_R$}\rightarrow
\fcat{cat$^{op}(J)$}{mod$_R$}$ acts by restriction, and admits a
right adjoint $G_*$, known as the {\it right Kan extension}
\cite{ma:cwm}. In particular, $G_*F^J$ is given on $\sigma\in K$ by
$\lim F^{\partial(\mu)}$ when $\sigma=\mu$, and $F^\sigma$
otherwise.

But $F^\mu\rightarrow\lim F^{\partial(\mu)}$ is an epimorphism, by
Lemma \ref{dkfat}, so the natural transformation $F^K\rightarrow
G_*F^J$ is epic on every face of $K$, and its kernel $H$ is zero on
every face except $\mu$. We acquire a short exact sequence of functors
\[
0\longrightarrow H\longrightarrow F^K\longrightarrow
G_*F^J\longrightarrow 0,
\]
which induces a long exact sequence of higher limits. By Lemma
\ref{atomic}, this collapses to a sequence of isomorphisms
\begin{equation}\label{isos}
{\textstyle\lim^i}F^K\;\cong\;{\textstyle\lim^i}G_*F^J,
\end{equation}
for $i\geq 1$. We now apply the composition of functors spectral
sequence \cite{caei:ha}, \cite{gazi:cfh}
\[
{\textstyle\lim^n}G_*^iF^J\Longrightarrow{\textstyle\lim^{n+i}} F^J.
\]
Here $G_*^i$ denotes the $i$th derived functor of $G_*$; it may be
evaluated on any face $\sigma$ of $K$ as $\lim^iF^{\partial(\sigma)}$,
and therefore vanishes for $i>0$, by inductive hypothesis. So the
spectral sequence collapses onto the first row of the $E^2$ page, from
which we obtain isomorphisms $\lim^nG_*F^J\cong\lim^nF^J$ for all $n\geq
0$. Since the inductive hypothesis applies to $J$, we deduce that
$\lim^nG_*F^J=0$ for every $n>0$. Combining this with \eqref{isos}
concludes the proof.
\end{proof}
\begin{cor}\label{codesharp}
The \bkss\ for $B^K$ collapses at the $E_2$ page; it is
concentrated along the vertical axis, and given by
\[
{\textstyle\lim^i}H^j(B^K;R)=
\begin{cases}
\lim H^j(B^K;R)&\text{if $i=0$}\\ 0&\text{otherwise}.
\end{cases}
\]
\end{cor}
\begin{proof}
The result follows immediately from Theorem \ref{hilimze} by letting
$D^K$ be $H^j(-;R)$ for every $j\geq 0$.
\end{proof}
Corollary \ref{codesharp} confirms that the edge homomorphism $h$ is
an isomorphism in $\cat{gca}_R$. When combined with \eqref{hrprojs}
and Proposition \ref{epi} it implies that the natural map
\[
R[K]\;\cong\;H^*(c(K);R)\stackrel{h}{\lra}\lim H^*(B^K;R),
\]
which is induced by the projections $R[K]\rightarrow S_R(\sigma)$,
is also an isomorphism. This may be proven directly, by refining the
methods of Proposition \ref{epi}. 

\section{Integral formality}\label{info}

In this section we study the formality of $c(K)$ over our arbitrary
commutative ring $R$, and construct a zig-zag of weak equivalences 
between the singular cochain algebra $C^*(c(K);R)$ and its cohomology
ring. 

We work in the model category $\cat{dga}_R$ of differential graded
$R$-algebras, whose differentials have degree $+1$; morphisms which
induces isomorphisms in cohomology are known as {\it
quasi-isomorphisms}. The model structure arises by interpreting
$\cat{dga}_R$ as the category of monoids in the monoidal model
category $\cat{dgmod}_R$ of unbounded cochain complexes over $R$. The
latter is isomorphic to Hovey's category of unbounded chain complexes
\cite{ho:mc}, and the model structure is induced on $\cat{dga}_R$ by
checking that it satisfies the monoid axiom of \cite{scsh:amm}. As
Schwede and Shipley confirm, weak equivalences are quasi-isomorphisms
and fibrations are epimorphisms. Cofibrations are defined by the
appropriate lifting property, and are necessarily degreewise split
injections. We emphasise that the objects of $\cat{dga}_R$ need not be
commutative.

A differential graded $R$-algebra $C^*$ is {\it formal in
$\cat{dga}_R$} whenever there is a zig-zag of quasi-isomorphisms
\begin{equation}\label{fodga}
H(C^*)\stackrel{\sim}{\longrightarrow}\dots
\stackrel{\sim}{\longleftarrow}C^*
\end{equation}
in $\cat{dga}_R$, where we follow the standard convention of assigning
the zero differential to the cohomology algebra $H(C^*)$. Our aim is to
show that the cochain algebra $C^*(c(K);R)$ is always formal in
$\cat{dga}_R$. This extends Franz's result \cite{fr:ics}, which only
applies to complexes arising from smooth fans.

We begin by choosing $D$ to be the $j$-dimensional cochain functor
$C^j(-;R)$ in \eqref{dkdk}, thus creating twin diagrams
$\big(C^j(B_K;R),C^j(B^K;R)\big)$ for each $j\geq 0$. As in
\eqref{grcobk}, we may consider the graded version $C^*(-;R)$ in
$\cat{dgmod}_R$. In fact its values are always $R$-algebras, with
respect to the cup product of cochains. The product is not commutative,
but the procedure for forming the limit of a $\cat{dga}_R$-diagram
remains the same; work in $\cat{dgmod}_R$, and superimpose the induced
multiplicative structure.

For the \EM\ space $B$, we let $v$ denote a generator of $H^2(B;R)$,
which is isomorphic to $R$. We choose a cocycle $\psi_v$ representing
$v$ in $C^2(B;R)$, and define a homomorphism $\psi\colon
H^*(B;R)\rightarrow C^*(B;R)$ by $\psi(v^k)=(\psi_v)^k$, for all
$k\geq 0$. By construction, $\psi$ is multiplicative, and is a
quasi-isomorphism in $\cat{dga}_R$. 
In order to extend this procedure we introduce a quasi-isomorphism
$\kappa$, defined by composition with the K\"unneth isomorphism as
\begin{equation}\label{comp}
H^*(B^V;R)\stackrel{\cong}{\llongrightarrow}H^*(B;R)^{\otimes V}
\stackrel{\psi^{\otimes}}{\llongrightarrow}C^*(B;R)^{\otimes V};
\end{equation} 
$\kappa$ is also multiplicative. There is a further zig-zag of
quasi-isomorphims
\begin{equation}\label{zzz}
\Hom(C_*(B);R)^{\otimes V}\longrightarrow
\Hom(C_*(B)^{\otimes V};R)\stackrel{ez}{\llongleftarrow}
\Hom(C_*(B^V);R),
\end{equation}
in which {\it ez\/} is the dual of Eilenberg-Zilber map. Both arrows
lie in $\cat{dga}_R$, so we may combine \eqref{comp} and \eqref{zzz}
to create the zig-zag
\begin{equation}\label{zz1}
H^*(B^V;R)\stackrel{\kappa}{\llongrightarrow}C^*(B;R)^{\otimes V}
\longrightarrow\Hom(C_*(B)^{\otimes
V};R)\stackrel{ez}{\llongleftarrow} C^*(B^V;R).
\end{equation}
This confirms that $C^*(B^n;R)$ is formal in $\cat{dga}_R$ for every
$n\geq 1$.

For each $B^\sigma$, we may project \eqref{zz1} onto the
corresponding zig-zag of quasi-isomorphisms. The results are
compatible by naturality, and so provide morphisms
\[
H^*(B^K;R)\stackrel{\kappa}{\llongrightarrow}C^*(B;R)^{\otimes K}
\longrightarrow\Hom(C_*(B)^{\otimes K};R)
\stackrel{ez}\llongleftarrow C^*(B^K;R)
\]
of $\cat{cat}^{op}(K)$-diagrams. Taking limits in $\cat{dga}_R$
yields the zig-zag
\begin{equation}\label{bigzz}
\begin{split}
\lim H^*(B^K;R)\stackrel{\kappa}{\llongrightarrow}\lim
C^*&(B;R)^{\otimes K}\longrightarrow\\
\lim\big(\Hom&(C_*(B)^{\otimes K};R)\big) 
\stackrel{ez}\llongleftarrow \lim C^*(B^K;R).
\end{split}
\end{equation}

\begin{lem}\label{limquiso}
All three homomorphisms of \eqref{bigzz} are quasi-isomorphisms in
$\cat{dga}_R$.
\end{lem}
\begin{proof}
A diagram $D^K\colon\cat{cat}^{op}(K)\rightarrow\cat{dga}_R$ is
fibrant whenever the projection onto the constant diagram $0$ is a
fibration. By \eqref{diagf} this occurs precisely when $D^K$ is fat,
and therefore holds for $H^*(B^K;R)$ and $C^*(B^K;R)$ by Lemma
\ref{dkfat}; it follows for $H^*(B;R)^{\otimes K}$ by the K\"unneth
isomorphism. So far as $C^*(B;R)^{\otimes K}$ is concerned, we note
that singular cochains determine a pair of twin diagrams $(C_{\otimes
K},C^{\otimes K})$ in $\cat{dga}_R$. Both functors assign
$C^*(B;R)^{\otimes\sigma}$ to the face $\sigma$. The value of
$C_{\otimes K}$ on $i_{\sigma,\tau}$ is the inclusion
$C^*(B;R)^{\otimes\sigma}\rightarrow C^*(B;R)^{\otimes\tau}$, and the
value of $C^{\otimes K}$ on $p_{\tau,\sigma}$ is the projection
$C^*(B;R)^{\otimes\tau}\rightarrow C^*(B;R)^{\otimes\sigma}$; the
latter requires the augmentation induced by the base point of
$B$. Hence $C^{\otimes K}$ is also fat, and $C^*(B;R)^{\otimes K}$ is
fibrant. Similar remarks apply to $\Hom(C_*(B)^{\otimes K};R)$.

The homomorphisms in question are therefore objectwise equivalences of
fibrant diagrams, and induces weak equivalences of limits by
\cite{ho:mc}.
\end{proof}

\begin{lem}\label{cocoli}
The natural homomorphism $g\colon C^*(c(K);R)\rightarrow\lim
C^*(B^K;R)$ is a quasi-isomorphism in $\cat{dga}_R$.
\end{lem}
\begin{proof} The edge isomorphism $h$ of \eqref{hrprojs} and Corollary  
\ref{codesharp} factorises as
\[ 
H(C^*(c(K);R))\stackrel{H(g)}{\llongrightarrow}H(\lim C^*(B^K;R))
\stackrel{l}{\longrightarrow}\lim H^*(B^K;R) 
\] 
in $\cat{dga}_R$, where $l$ is induced by the collection of compatible
homomorphisms 
\[ 
H(\lim C^*(B^K;R))\longrightarrow H^*(B^\sigma;R).  
\] 
Now let $d$ be the differential on $C^j(-;R)$ for every $j\geq 0$, and
define the cycle and boundary functors
$Z^j,\,I^j\colon\cat{top}\rightarrow \cat{mod}_R$ as the kernel and
image of $d$ respectively. They determine twin diagrams, and therefore
fat functors
$Z^K,\,I^K\colon\cat{cat}^{op}(K)\rightarrow\cat{mod}_R$. Theorem
\ref{hilimze} confirms that $\lim^iZ^j(B^K;R)=\lim^iI^j(B^K;R)=0$ for
all $i>0$ and $j\geq 0$.  It follows immediately that $l$ is an
isomorphism, and therefore that $H(g)$ is an isomorphism, as sought.
\end{proof}

We may now complete our analysis of $C^*(c(K);R)$.
\begin{thm}\label{cckrfo}
The differential graded $R$-algebra $C^*(c(K);R)$ is formal in
$\cat{dga}_R$.
\end{thm}
\begin{proof}
Combining Corollary \ref{codesharp} with Lemmas \ref{limquiso} and
\ref{cocoli} yields a zig-zag
\[
H^*(c(K);R)\stackrel{h}{\longrightarrow}\lim
H^*(B^K;R)\longrightarrow \dots\longleftarrow\lim
C^*(B^K;R)\stackrel{g}{\longleftarrow}C^*(c(K);R)
\]
of quasi-isomorphisms, as required by \eqref{fodga}.
\end{proof}

\begin{rem}\label{genintfor}
The proof of Theorem \ref{cckrfo} extends to exponential diagrams
$X^K$ for which $C^*(X;R)$ is formal in $\cat{dga}_R$ and the
K\"unneth isomorphism $H^*(X^V;R)\cong H^*(X;R)^{\otimes V}$ holds. We
replace $\psi$ in \eqref{zz1} by the corresponding zig-zag
$H^*(X;R)\stackrel{\sim}{\longrightarrow}\dots
\stackrel{\sim}{\longleftarrow}C^*(X;R)$ of quasi-isomorphisms, and
repeat the remainder of the argument above.
\end{rem}

\section{Rational formality}\label{rafo}

In our final section we turn to the rational case $R=\bQ$, and confirm
the formality of Sullivan's algebra of rational cochains on $c(K)$ in
the commutative setting. This involves stricter conditions than those
for general $R$, and has deeper topological consequences. In particular,
it leads us to a minimal model whenever $\bQ[K]$ is a complete
intersection ring, and thence to rational uniqueness. We refer
readers to Bousfield and Gugenheim \cite{bogu:pdt} for details of the
model category of differential graded commutative $\bQ$-algebras, and to
F\'elix, Halperin and Thomas \cite{fehath:rht} for background
information on rational homotopy theory.

We begin by replacing $C^*(X;R)$ with Sullivan's rational algebra
$\APL(X)$ of polynomial forms\cite{fehath:rht}. The commutativity of
the latter is crucial, and suggests we work in the category
$\cat{dgca}_\bQ$ of differential graded commutative $\bQ$-algebras
\cite{bogu:pdt}. The existence of a model structure is assured by
working over a field; as before, weak equivalences are
quasi-isomorphisms, fibrations are epimorphisms, and cofibrations are
defined by the appropriate lifting properties.

For each $s\geq 0$, we write the differential algebra of rational
polynomial forms on the standard $s$-simplex as $\nabla_s(*)$. It is an
object of $\cat{dgca}_\bQ$.  For each $t\geq 0$, the forms of
dimension $t$ define a simplicial vector space $\nabla_\bullet(t)$ over
$\bQ$, and $\nabla_\bullet(*)$ becomes a simplicial object in
$\cat{dgca}_\bQ$. So
\[
A^*(Y_\bullet)\letbe\cat{sset}(Y_\bullet,\nabla_\bullet(*))
\]
is also an object of $\cat{dgca}_\bQ$, which is weakly equivalent to
the normalised cochain complex $N^*(Y_\bullet;\bQ)$. Then $\APL(X)$ is
defined as $A^*(S_\bullet X)$, where $S_\bullet$ denotes the total
singular complex functor $\cat{top}\rightarrow\cat{sset}$. The PL de
Rham Theorem \cite{bogu:pdt} asserts that the cohomology algebra
$H(\APL(X))$ is naturally isomorphic to $H^*(X;\bQ)$. As usual, we
invest cohomology algebras with the zero differential.

A differential graded commutative $\bQ$-algebras $A^*$ is {\it formal in
$\cat{dgca}_\bQ$} whenever there is a zig-zag of quasi-isomorphisms
\begin{equation}\label{foquis}
H(A^*)\stackrel{\sim}{\longrightarrow}\dots
\stackrel{\sim}{\longleftarrow}A^*
\end{equation}
in $\cat{dgca}_\bQ$. A topological space $X$ is {\it rationally
formal\/} whenever $\APL(X)$ satisfies \eqref{foquis}.

For any such $X$, a minimal Sullivan model may be constructed directly
from the algebra $H^*(X;\bQ)$. Our remaining goal is therefore to show
that $c(K)$ is rationally formal, and to consider the implications for
the uniqueness of spaces $X$ realising $\bQ[K]$. The proof parallels
that for general $R$, but the need to respect commutativity forces
several changes of detail.

We choose $D$ to be $\APL$ in \eqref{dkdk}, creating twin diagrams
$(\APL(B_K),\APL(B^K))$ in $\cat{dgca}_\bQ$. As before, we form
limits by working in $\cat{dgmod}_Q$, and superimposing the induced
multiplicative structure. Applying cohomology yields the twins
$(H^*(B_K;\bQ),H^*(B^K;\bQ))$, whose value on each face $\sigma$ is
$S_\bQ(\sigma)$ in $\cat{dgmod}_Q$. Both $\APL(B^K)$ and
$H^*(B^K;\bQ)$ are fat, by Lemma \ref{dkfat}.

Using the fact that $H(\APL(B^V))$ is isomorphic to $S_\bQ(V)$, we
choose cocycles $\phi_j$ in $\APL(B^V)$ representing $v_j$ for every
$1\leq j\leq m$. We may then define a homomorphism
\begin{equation}\label{defphi}
\phi\colon H^*(B^V;\bQ)\longrightarrow\APL(B^V)
\end{equation}
by $\phi(v_j)=\phi_j$, because $\APL(B^V)$ is commutative. Moreover,
$\phi$ is a quasi-isomorph\-ism, reflecting the rational formality of the
\EM\ space $H(\bZ^V;2)$. By restriction, we interpret the $\phi_j$ as
cocycles in $\APL(B^\sigma)$ for every face $\sigma$. They then
represent $v_j$ in $S_\bQ(\sigma)$ when $\sigma$ contains $v_j$, and $0$
otherwise. We obtain compatible quasi-isomorphisms on each $B^\sigma$,
which combine to create a map  
\[
\phi\colon H^*(B^K;\bQ)\longrightarrow\APL(B^K)
\]
of $\cat{cat}^{op}(K)$-diagrams in $\cat{dgca}_\bQ$. It is an
objectwise weak equivalence. Taking limits yields a homomorphism
\[
l(\phi)\colon\lim H^*(B^K;\bQ)\longrightarrow\lim\APL(B^K)
\]
of differential graded commutative algebras over $\bQ$.
\begin{lem}\label{phiso}
The homomorphism $l(\phi)$ is a quasi-isomorphism in
$\cat{dgca}_\bQ$.
\end{lem}
\begin{proof}
Both diagrams are fat, and therefore fibrant by \eqref{diagf}. So
$\phi$ induces a weak equivalence of limits.
\end{proof}

Because the $B^\sigma$ are \EM\ spaces, it is convenient to complete our
proof of Theorem \ref{xkformal} in terms of simplicial sets. We may
then take advantage of the fact that $A^*$ converts colimits in
$\cat{sset}$ to limits in $\cat{dgca}_\bQ$, for reasons which are
purely set-theoretic.

We denote the realisation functor $\cat{sset}\rightarrow\cat{top}$ by
$|\!-\!|$. Given an arbitrary simplicial set $Y_\bullet$, there is a
quasi-isomorphism
\begin{equation}\label{aplsset}
\APL(|Y_\bullet|)\,=\, A^*(S_\bullet|Y_\bullet|)
\stackrel{\sim}{\lra}A^*(Y_\bullet),
\end{equation}
induced by the natural equivalence $Y_\bullet\rightarrow
S_\bullet|Y_\bullet|$. For each face $\sigma$ of $K$, we choose
$|H_\bullet(\bZ^\sigma;2)|$ as our model for $B^\sigma$; it is
well-pointed, by the cofibration induced by the inclusion of the
trivial subgroup $\{0\}\rightarrow\bZ^\sigma$. We write $H_\bullet^K$
for the corresponding diagram of simplicial sets, which takes the
value $H_\bullet(\bZ^\sigma;2)$ on $\sigma$.
\begin{thm}\label{xkformal}
The space $c(K)$ is rationally formal.
\end{thm}
\begin{proof}
Since realisation is left adjoint to $S_\bullet$, it commutes with
colimits. So we may write
\[
c(K)=\col|H_\bullet^K|\cong|\col H_\bullet^K|,
\]
where the second colimit is taken in $\cat{sset}$. Applying
\eqref{aplsset} gives a zig-zag
\begin{equation}\label{zz2}
\lim\APL(B^K)\stackrel{\sim}{\longrightarrow}\lim A^*(H_\bullet^K)
\stackrel{\cong}{\longrightarrow}A^*(\col H_\bullet^K)
\stackrel{\sim}{\longleftarrow}\APL(c(K)),
\end{equation}
where the central isomorphism follows from the property of $A^*$
described above. Combining \eqref{zz2} with Corollary \ref{codesharp}
and Lemma \ref{phiso} yields a zig-zag 
\[
H^*(c(K);\bQ)\stackrel{h}{\longrightarrow}\lim
H^*(B^K;\bQ)\stackrel{l(\phi)}{\llongrightarrow}\lim\APL(B^K)
\stackrel{\sim}{\longrightarrow}\dots
\stackrel{\sim}{\longleftarrow}\APL(c(K))
\]
of quasi-isomorphisms in $\cat{dgca}_\bQ$. The result follows from
\eqref{foquis}.
\end{proof}

\begin{rem}
By analogy with Remark \ref{genintfor}, the proof of Theorem
\ref{xkformal} extends to exponential diagrams $X^K$ for which $X$ is
rationally formal; the K\"unneth isomorphism holds automatically,
because we are working over $\bQ$. The product $\APL(X)^{\otimes
V}\rightarrow\APL(X^V)$ of the maps induced by projection is a
quasi-isomorphism, and is natural with respect to projection and
inclusion of coordinates. So we may replace $\phi$ in \eqref{defphi}
by the corresponding zig-zag of quasi-isomorphisms
\[
H^*(X^V;\bQ)\stackrel{\cong}{\lra}H^*(X;\bQ)^{\otimes V}
\stackrel{\sim}{\lra}\dots\stackrel{\sim}{\longleftarrow}
\APL(X)^{\otimes V}\stackrel{\sim}{\lra}\APL(X^V), 
\] 
and proceed with the remainder of the argument above.
\end{rem}

Theorem \ref{xkformal} confirms that a minimal Sullivan model for
$c(K)$ may be constructed directly from $\bQ[K]$. It consists of an
acyclic fibration
\begin{equation}\label{minfib}
\eta\colon(S_\bQ(W(K)),d)\longrightarrow H^*(c(K);\bQ),
\end{equation}
where $W(K)$ is an appropriately graded set of generators
(necessarily exterior in odd dimensions), and provides a cofibrant
replacement for $H^*(c(K);\bQ)$ in $\cat{dgca}_\bQ$. In general,
$W(K)$ is not easy to describe, although special cases such as
Example \ref{minmods} below are straightforward, and lead to our
uniqueness result. The properties of $W(K)$ are linked to those of
the loop space $\varOmega c(K)$, whose study was begun in
\cite{paravo:csa}; we expect to return to this relationship in
future.

The principal calculational tool of rational homotopy theory is the
Sullivan-de Rham equivalence of homotopy categories, which asserts that
Bousfield and Gugenheim's adjoint pair of derived functors
\[
ho\,\cat{sset}_\bQ\que ho\,\cat{dgca}_\bQ
\]
restrict to inverse equivalences between certain full subcategories
\cite{bogu:pdt}. These are given by nilpotent simplicial sets of
finite type, and algebras which are equivalent to minimal algebras
with finitely many generators in each dimension.

In particular, the equivalence identifies homotopy classes of maps
$[c(K)_0,c(L)_0]$ with homotopy classes of morphisms $[S_\bQ
(W(L)),S_\bQ(W(K))]$. Since every object of $\cat{dgca}_\bQ$ is
fibrant, it suffices to consider homotopy classes of the form
$[S_\bQ(W(L)),H^*(c(K);\bQ)]$; of course $S_\bQ(W(L))$ cannot be
substituted similarly, because $H^*(c(L);\bQ)$ is not usually
cofibrant. Nevertheless, the function
\begin{equation}\label{takeco}
[S_\bQ(W(L)),H^*(c(K);\bQ)]\longrightarrow
\cat{gca}_\bQ(H^*(c(L);\bQ),H^*(c(K);\bQ))
\end{equation}
induced by taking cohomology is always surjective, and it would be
of interest to understand its kernel.

\begin{exa}\label{minmods}
Let $(\lambda(k):1\leq k\leq t)$ be a sequence of disjoint subsets
of $V$, where $\lambda(k)$ has cardinality $n(k)$, and define $L$ to
be the subcomplex of $\varDelta(V)$ obtained by deleting all faces
containing one or more of the $\lambda(k)$. We write {\em
$\lambda\letbe\cup_k\lambda(k)$}, and $|\lambda|=n$. Over any
commutative ring $R$, the Stanley-Reisner algebra $R[L]$ is given by
quotienting out the regular sequence $v_{\lambda(1)}$, \dots
$v_{\lambda(t)}$ from $S_R(V)$.

Over $\bQ$, the generating set $W(L)$ of \eqref{minfib} consists of $V$ 
in dimension $2$, and elements $w(k)$ in dimension $2n(k)-1$, for $1\leq
k\leq t$; the differential is given by $dv_j=0$ for all $j$, and
$dw(k)=v_{\lambda(k)}$. The fibration $\eta$ identifies the vertices
$V$ in dimension $2$, and annihilates every $w(k)$.
\end{exa}

Since the elements $w(k)$ are odd dimensional, every
$\cat{dgca}_\bQ$-morphism 
\[
S_\bQ(W(L))\rightarrow H^*(c(K);\bQ)
\]
is determined by its effect on $V$, and the function \eqref{takeco} is
bijective. It follows that $\HA(c(L)_0)$ (as defined in Section
\ref{intro}) is isomorphic to the group of algebra automorphisms
$\Aut(\bQ[L])$, and is therefore a subgroup of $\GL(m,\bQ)$. It
contains all matrices of the form $\left(\smallmatrix
M&0\\N&\varSigma\endsmallmatrix\right)$, where $M\in\GL(m-n,\bQ)$ acts
on $\bQ^{V\setminus\lambda}$, and $\varSigma$ permutes the elements of
$\lambda$. The permutations act on the elements of each individual
$\lambda(k)$, and interchange those $\lambda(k)$ which are of common
cardinality. We conjecture that every automorphism of $\bQ[L]$ may be
represented by such a matrix.

It is convenient to refer to $L$ as a {\it complete intersection
complex\/} whenever $\bQ[L]$ is a complete intersection ring. A
straightforward application of \cite[Theorem 2.3.3]{brhe:cmr} shows
that this occurs if and only if $L$ takes the form of Example
\ref{minmods}.

Our final result concerns uniqueness, and is a consequence of the fact
that complete interseection complexes are a special case of Sullivan's
examples \cite[page 317]{su:ict}. Following Sullivan, it may be
summarised by the statement that $c(L)$ and $\bQ[L]$ are {\it
intrinsically formal\/} for any such $L$.
\begin{prop}\label{unique}
Let $X$ be a nilpotent CW-complex of finite type, which realises a
complete intersection ring $\bQ[L]$; then the rationalisations of $X$
and $c(L)$ are homotopy equivalent.
\end{prop}
\begin{proof}
We write the generators of the cohomology algebra $H^*(X;\bQ)$ as
$v'_j$, where $1\leq j\leq m$, and choose representing cocycles $\phi_j$
in $\APL(X)$. Each monomial $v'_M$ is therefore represented by the
corresponding product $\phi_M$. It follows that $\phi_{\lambda(k)}$ is a
coboundary, and there exist elements $\theta(k)$ such that
$d\theta(k)=\phi_{\lambda(k)}$ in $\APL(X)$ for all $1\leq k\leq t$.

We define a $\cat{dgca}_\bQ$-morphism $\eta\colon
S_\bQ(W(L))\rightarrow\APL(X)$ by setting $\eta(v_j)\letbe\phi_j$ and
$\eta(w(k))\letbe\theta(k)$, for $1\leq j\leq m$ and $1\leq k\leq t$
resepectively. So $\eta$ is a weak equivalence, and is a minimal model
for $X$. Hence $X$ and $c(L)$ have isomorphic minimal models, and the
result follows.
\end{proof}

\Addresses\recd

\end{document}

%% file: agtout.tex
\def\ifplaintex{\expandafter\ifx\csname documentclass\endcsname\relax}

\def\gtp{{\mathsurround=0pt\it $\cal G\mskip-2mu$eometry \&\ 
$\cal T\!\!$opology $\cal P\!$ublications}}  

\def\recd{{\small Received:\qua\receiveddate\ifx\reviseddate\relax
\else\qquad Revised:\qua\reviseddate\fi\par}} 

\def\lognumber#1{\def\thelognumber{#1}}
\def\volumenumber#1{\def\thevolumenumber{#1}}
\def\volumeyear#1{\def\thevolumeyear{#1}}
\def\papernumber#1{\def\thepapernumber{#1}}
\def\pagenumbers#1#2{\def\startpage{#1}\def\finishpage{#2}}
\def\published#1{\def\publishdate{#1}}

\def\received#1{\def\receiveddate{#1}}
\def\revised#1{\def\reviseddate{#1}}
\def\accepted#1{\def\accepteddate{#1}}

\def\asciiaddress#1{\def\theasciiaddress{#1}}
\def\asciiemail#1{\def\theasciiemail{#1}}

\long\def\asciiabstract#1{\long\def\theasciiabstract{#1}}
\def\asciikeywords#1{\def\theasciikeywords{#1}}

\let\\\par\let\thelognumber\relax\let\thevolumenumber\relax
\let\thepapernumber\relax\let\thevolumeyear\relax\let\startpage\relax
\let\finishpage\relax\let\publishdate\relax\let\receiveddate\relax
\let\reviseddate\relax\let\accepteddate\relax\let\theasciititle\relax
\let\theasciiauthors\relax\let\theasciiaddress\relax
\let\theasciiabstract\relax\let\theasciikeywords\relax

\let\theasciiemail\relax

\ifplaintex
\font\logobig=cmssbx10 scaled 3836
\font\logomed=cmssbx10 scaled 2557
\else
\font\logobig=cmssbx10 scaled 4200
\font\logomed=cmssbx10 scaled 2800
\fi

\long\def\makeagttitle{   
\count0=\startpage
\agt\hfill      
\hbox to 45truept{\vbox to 0pt{\vglue -13truept{\logomed A\kern -.37em{\logobig 
T}\kern -.38em G}\vss}\hss}
\break
{\small Volume \thevolumenumber\ (\thevolumeyear)
\startpage--\finishpage\nl
Published: \publishdate}

\vglue .25truein

{\parskip=0pt\leftskip 0pt plus
1fil\def\\{\par\smallskip}{\Large\bf\thetitle}\par\medskip} \vglue
0.05truein

{\parskip=0pt\leftskip 0pt plus 1fil\def\\{\par}{\sc\theauthors}
\par\medskip}%
 
\vglue 0.03truein 

{\small\leftskip 25truept\rightskip 25truept{\bf Abstract}\stdspace\theabstract

{\bf AMS Classification}\stdspace\theprimaryclass
\ifx\thesecondaryclass\relax\else; \thesecondaryclass\fi\par
{\bf Keywords}\stdspace \thekeywords\par}\vglue 7truept

}   

\ifplaintex
\font\phead=cmsl9 scaled 950
\font\pnum=cmbx10 scaled 913
\font\pfoot=cmsl9 scaled 950
\headline{\vbox to 0pt{\vskip -4.5mm\line{\small\phead\ifnum
\count0=\startpage ISSN 1472-2739 (on-line) 1472-2747 (printed)
\hfill {\pnum\folio}\else\ifodd\count0\def\\{ }%
\ifx\theshorttitle\relax\thetitle\else\theshorttitle\fi\hfill{\pnum\folio}
\else\def\\{ and }{\pnum\folio}\hfill\ifx\theshortauthors\relax\theauthors
\else\theshortauthors\fi\fi\fi}\vss}}
\footline{\vbox to 0pt{\vglue 0mm\line{\small\pfoot\ifnum\count0=\startpage
\copyright\ \gtp\hfill\else
\agt, Volume \thevolumenumber\ (\thevolumeyear)\hfill\fi}\vss}}
\else
\font=cmsl9 scaled 1050
\font\lnum=cmbx10 
\font=cmsl9 scaled 1050
\makeatletter
\def\@oddhead{{\small\lhead\ifnum\count0=\startpage ISSN 1472-2739 
(on-line) 1472-2747 (printed)\hfill {\lnum\number\count0}\else\ifodd\count0
\def\\{ }\ifx\theshorttitle\relax \thetitle \else\theshorttitle\fi\hfill
{\lnum\number\count0}\else\def\\{ and }{\lnum\number\count0}
\hfill\ifx\theshortauthors\relax 
\theauthors\else\theshortauthors\fi\fi\fi}}\def\@evenhead{\@oddhead}
\def\@oddfoot{\small\lfoot\ifnum\count0=\startpage\copyright\ \gtp\hfill\else
\agt, Volume \thevolumenumber\ (\thevolumeyear)\hfill\fi}
\def\@evenfoot{\@oddfoot}
\makeatother
\fi
\let\maketitlepage\makeagttitle
\let\makeshorttitle\maketitlepage
\let\maketitle\maketitlepage

\newwrite\gtoutfile
\long\gdef\makeheadfile{  
{\def\\{, }\def\s{ }
\immediate\openout\gtoutfile head.xxx
\immediate\write\gtoutfile{Proxy-for: \ifx\theasciiauthors\relax
\theauthors\else\theasciiauthors\fi\s<\ifx\theasciiemail\relax\theemail\else\theasciiemail\fi>}
\immediate\write\gtoutfile{\noexpand\\}
\immediate\write\gtoutfile{Authors: \ifx\theasciiauthors\relax
\theauthors\else\theasciiauthors\fi}
{\def\\{ }\immediate\write\gtoutfile{Title: \ifx\theasciititle\relax
\thetitle\else\theasciititle\fi}}
\immediate\write\gtoutfile{Subj-class: GT or SG, GR etc}
\immediate\write\gtoutfile{MSC-class: \theprimaryclass\ifx\thesecondaryclass\relax\else, \thesecondaryclass\fi}
\immediate\write\gtoutfile{Journal-ref: Algebr. Geom. Topol. \thevolumenumber\s
(\thevolumeyear) \startpage-\finishpage}
\immediate\write\gtoutfile{Comments: Published by Algebraic and
Geometric Topology at}
\immediate\write\gtoutfile{\s\s\s  http://www.maths.warwick.ac.uk/agt/AGTVol\thevolumenumber/agt-\thevolumenumber-\thepapernumber.abs.html}
\immediate\write\gtoutfile{\noexpand\\}
\immediate\write\gtoutfile{}
\ifx\theasciiabstract\relax
\immediate\write\gtoutfile{\theabstract}\else
\immediate\write\gtoutfile{\theasciiabstract}\fi
\immediate\write\gtoutfile{}
\immediate\write\gtoutfile{\noexpand\\}
\immediate\write\gtoutfile{}
\immediate\closeout\gtoutfile}}  

\def\maketitlepage{\makeagttitle\makeheadfile}
\let\makeshorttitle\maketitlepage
\let\maketitle\maketitlepage